\documentclass[11pt]{amsart}

\usepackage{amssymb}
\usepackage{amsmath, amsthm, amsopn, amsfonts}
\usepackage[dvips]{graphicx}
\usepackage{graphpap}
\usepackage[francais, english]{babel}

\def\R{\mathbb R}

\def\N{\mathbb N}

\newtheorem{proposition}{Proposition}

\newtheorem{remarque}[proposition]{Remark}

\newtheorem{theo}{Theorem}[section]
\newtheorem{rem}{Remark}[section]
\newtheorem{lem}{Lemma}[section]

\numberwithin{equation}{section}
\numberwithin{proposition}{section}
\newcommand{\rien}[1]{\relax}% DIVERS

\begin{document}
\selectlanguage{english}

\title
[ On KP type equations]{\bf Remarks on the mass constraint for KP type equations }
\author{L.~Molinet} 
\address{L.A.G.A., Institut Galil\'ee, Universit\'e Paris 13, 93430 Villetaneuse, France }
\author{J. C.~Saut}
\address{ Universit\'e de Paris-Sud, UMR de Math\'ematiques, B\^at. 425, 91405 Orsay Cedex, France }
\author {N.~Tzvetkov}
\address
{D\'epartement de Math\'ematiques, Universit\'e Lille I, 59 655 Villeneuve d'Ascq Cedex, France}
%\date{}
\begin{abstract}
For a rather general class of equations of Kadomtsev-Petviashvili (KP) type, we prove that the zero-mass (in $x$)
constraint is satisfied at any non zero time even if it is not
satisfied at initial time zero.
Our results are based on a precise analysis of the fundamental solution of the
linear part and its anti $x$-derivative.
\end{abstract} 
\maketitle
\section{Introduction}
Kadomtsev-Petviashvili equations are universal for the modelling of the propagation of long weakly dispersive waves which propagate
essentially in one direction with weak transverse effects. As explained in the pioneering paper of Kadomtsev and Petviashvili 
\cite{KP}, they are (formally) obtained in the following way. 
We start from a one dimensional long-wave dispersive equation which is of Korteweg- de Vries (KdV) type, i.e
\begin{equation}\label{1.1}
u_{t}+u_{x}+uu_{x}-Lu_{x}=0,\quad u=u(t,x),\,x\in\R,\, t\geq 0\, .
\end{equation}
In (\ref{1.1}) $L$ is a (possibly nonlocal) operator, defined in Fourier variable by
\begin{equation}\label{1.3}
\widehat{Lf}(\xi)=c(\xi)\widehat{f}(\xi),
\end{equation}
where $c$ is a real function which is the phase velocity.
For instance, the case $c(\xi)=\pm\xi^2$ ($L=\mp\partial_{x}^{2}$) corresponds to the classical KdV equation.
In the context of water waves, the sign of $c(\xi)$ depends on the surface tension parameter. 
The case $c(\xi)=|\xi|$ $(L=H\partial_{x})$ corresponds to the Benjamin-Ono equation, etc ...

As observed in \cite{KP} the correction to (\ref{1.1}) due to weak transverse effects are independent of the
dispersion in $x$ and is related only to the finite propagation speed properties of the transport operator 
$M=\partial_{t}+\partial_{x}$. Recall that $M$ gives rise to one directional waves moving to the right with speed one,
i.e. a profile $\varphi(x)$ evolves under the flow of $M$ as $\varphi(x-t)$.
A weak transverse perturbation of $\varphi(x)$ 
is a two dimensional function $\psi(x,y)$ close to $\varphi(x)$, localized in the frequency region
$\Big|\frac{\eta}{\xi}\Big|\ll 1$, where $\xi$ and $\eta$ are the Fourier modes corresponding to $x$ and $y$ respectively.
We look for a two dimensional perturbation  $\widetilde{M}=\partial_{t}+\partial_{x}+\omega(D_x,D_y)$ of $M$ such that,
similarly to above, the profile of $\psi(x,y)$, does not change much
when evolving under the flow of $\widetilde{M}$. Here $\omega(D_x,D_y)$ denotes the Fourier multiplier with symbol
the real function $\omega(\xi,\eta)$. 
Natural generalizations of the flow of $M$ in $2D$ are the flows of the wave operators 
$\partial_{t}\pm\sqrt{-\Delta}$ which enjoy the finite propagation speed
property. 
Since
$$
\xi+\frac{1}{2}\xi^{-1}\eta^{2}\sim \pm\sqrt{\xi^2+\eta^2}, \quad  {\rm when }\quad\Big|\frac{\eta}{\xi}\Big|\ll 1\, .
$$
we deduce that
$$
\partial_{t}+\partial_{x}+\frac{1}{2}\partial_{x}^{-1}\partial_{y}^{2}\sim
\partial_{t}\pm \sqrt{-\Delta}
$$
which leads to the correction $\omega(D_x,D_y)=\frac{1}{2}\partial_{x}^{-1}\partial_{y}^{2}$. in (\ref{1.1}).

Of course when the transverse effect are two-dimensional, the correction is $\frac{1}{2}\partial_{x}^{-1}\Delta_{\perp}$,
where $\Delta_{\perp}=\partial_{y}^{2}+\partial_{z}^{2}$.

We are thus led to the following model in two dimensions 
\begin{equation}\label{1.4}
u_{t}+u_{x}+uu_{x}-Lu_{x}+\partial_{x}^{-1}\partial_{y}^{2}u=0,
\end{equation}
In equation (\ref{1.4}), it is implicitly assumed that the operator $\partial_{x}^{-1}\partial_{y}^{2}$
is well defined, which a priori imposes a constraint on the solution $u$, which, in some sense, has to be a $x$ derivative.
This is achieved for instance if  $u\in {\mathcal S}'(\R^2)$ is such that
\begin{equation}\label{1.6}
\xi_{1}^{-1}\,\xi_{2}^{2}\,\widehat{u}(t,\xi_1,\xi_2)\in {\mathcal S}'(\R^2)\, ,
\end{equation}
thus in particular if $\xi_{1}^{-1}\,\widehat{u}(t,\xi_1,\xi_2)\in {\mathcal S}'(\R^2)$.
Another possibility to fulfil the constraint is to write $u$ as
\begin{equation}\label{1.7}
u(t,x,y)=\frac{\partial}{\partial x}v(t,x,y)
\end{equation}
where $v$ is a continuous function having a classical derivative with respect
to $x$, which, for any fixed $y$ and $t\neq 0$ vanishes when $x\rightarrow \pm \infty$.
Thus one has
\begin{equation}\label{1.8}
\int_{-\infty}^{\infty}u(t,x,y)dx=0,\qquad y\in\R,\,\,\, t\neq 0
\end{equation}
in the sense of generalized Riemann integrals.
Of course the differentiated version of (\ref{1.4}), namely
\begin{equation}\label{1.9}
(u_{t}+u_{x}+uu_{x}-Lu_{x})_{x}+ \partial_{y}^{2}u=0
\end{equation}
can make sense without any constraint of type (\ref{1.6}) or  (\ref{1.8}) on $u$, and so does the Duhamel integral representation of (\ref{1.4}),
\begin{equation}\label{1.11}
u(t)=S(t)u_{0}-\int_{0}^{t}S(t-s)(u(s)u_{x}(s))ds
\end{equation}
where $S(t)$ denotes the (unitary in all Sobolev spaces $H^{s}(\R^2)$) group associated to (\ref{1.4}),
\begin{equation}\label{1.12}
S(t)=e^{-t(\partial_{x}-L\partial_{x}+\partial_{x}^{-1}\partial_{y}^{2})}\,\, .
\end{equation}

Let us notice, at this point, that models alternative to KdV type equations
(\ref{1.1}) are the equations of the Benjamin-Bona-Mahony (BBM) type
\begin{equation}\label{1.2}
u_{t}+u_{x}+uu_{x}+Lu_{t}=0
\end{equation}
with corresponding two dimensional ``KP-BBM type models'' (in the case $c(\xi)\geq 0$)
\begin{equation}\label{1.5}
u_{t}+u_{x}+uu_{x}+Lu_{t}+\partial_{x}^{-1}\partial_{y}^{2}u=0
\end{equation}
or
\begin{equation}\label{1.10}
(u_{t}+u_x+uu_{x}+Lu_{t})_{x}+ \partial_{y}^{2}u=0,
\end{equation}
and free group
$$
S(t)=e^{-t(I+L)^{-1}[\partial_{x}+ \partial_{x}^{-1}\partial_{y}^{2}]}\,\, .
$$
In view of the above discussion, all the results established for the Duhamel
form of KP type equations (e.g. those of Bourgain \cite{B}, Saut-Tzvetkov \cite{ST})
do not need any constraint on the initial data $u_0$.
It is then possible (see for instance \cite{MST1}) to check that the solution $u$ will
satisfy (\ref{1.9}) or (\ref{1.10}) in the distributional sense, but - not
a priori - the integrated forms (\ref{1.4}) or (\ref{1.5}).

On the other hand, a constraint has to be imposed when using the Hamiltonian
formulation of the equation. In fact, the Hamiltonian for (\ref{1.9}) is
\begin{equation}\label{1.13}
\frac{1}{2}\int\big[-u\,Lu+(\partial_{x}^{-1}u_y)^2+u^2+\frac{u^3}{3}\big]\, .
\end{equation} 
The Hamiltonian associated to (\ref{1.10}) is
\begin{equation}\label{1.14}
\frac{1}{2}\int\big[(\partial_{x}^{-1}u_y)^2+u^2+\frac{u^3}{3}\big]\, .
\end{equation} 
Therefore the global well-posedness results for KP-I obtained in
\cite{MST2},\cite{K} do need that the initial data satisfy (in particular) the
constraint $\partial_{x}^{-1}\partial_{y}u_0\in L^2(\R^2)$, and, this
constraint is preserved by the flow. Actually, the global results of
\cite{MST2,K} are partially based on higher order conservation laws of the
KP-I equation, and the constraint $\partial_{x}^{-2}\partial_{y}^2 u_0\in L^2(\R^2)$
is also needed (and preserved by the flow).

Actually, as noticed in \cite{MST2} there is a serious draw-back with the next
conservation laws of both KP-I and KP-II equations, starting with that
involving the $L^2$ norm of $\partial_{x}^{3}u$ which contains the $L^2$ norm of terms like
$(\partial_{x}^{-1}\partial_{y})u^2$.
For $u\in H^3(\R^2)$ this expression is meaningless in
$L^2(\R^2)$, unless $\partial_{y}(u^2)$ has zero mean in $x$,
which in turn implies that $u\equiv 0$.

The goal of this paper is to investigate the behaviour of a solution to
the general KP type equations (\ref{1.9}), (\ref{1.10}) which initially does
not satisfy the zero-mass constraint. We will show that in fact the zero-mass
constraint is satisfied at any non zero time $t$. 

At this point we should mention the papers \cite{AV,BPP1,BPP2} and especially
\cite{FS} where the inverse scattering transform machinery is used to solve
the Cauchy problem for KP-I and KP-II without the constraint. The more
complete and rigorous results are obtained in \cite{FS} (see also \cite{S}).
In the present work we consider a rather general class of KP or KP-BBM
equations and put the evidence on the key point which concerns only the linear
part : the fundamental solution of KP type equations is a $x$ derivative of a
$C^1$ with respect to $x$ continuous function which, for fixed $t\neq 0$ and $y$, tends to zero as
$x\rightarrow \pm \infty$. Thus its generalized Riemann integral in $x$
vanishes for all values of the transverse variable $y$ and of $t\neq 0$. A
similar property can be established for the solution of the non linear
problem. 

The paper is organized as follows. Section~2 deals with the KP type equations,
while Section~3 focuses on KP-BBM type equations. The last Section 4 reviews
briefly some extensions : the three dimensional case, non homogeneous
dispersive relations. 

In the sequel, different harmless numerical constants will be denoted
by $c$.
%%%%%%%%%%%%%%%%%%%%%%%%%%%%%%%%%%%%%%%%%%%%%%%%%%%%%%%%%%%%%%%%%%%%%%%%%%%%%%%%%%%%%%%%%%%%%%%%%%%%%%%%%%%%
%%%%%%%%%%%%%%%%%%%%%%%%%%%%%%%%%%%%%%%%%%%%%%%%%%%%%%%%%%%%%%%%%%%%%%%%%%%%%%%%%%%%%%%%%%%%%%%%%%%%%%%%%%%%
%%%%%%%%%%%%%%%%%%%%%%%%%%%%%%%%%%%%%%%%%%%%%%%%%%%%%%%%%%%%%%%%%%%%%%%%%%%%%%%%%%%%%%%%%%%%%%%%%%%%%%%%%%%%
%%%%%%%%%%%%%%%%%%%%%%%%%%%%%%%%%%%%%%%%%%%%%%%%%%%%%%%%%%%%%%%%%%%%%%%%%%%%%%%%%%%%%%%%%%%%%%%%%%%%%%%%%%%%
\section{KP type equations}
\subsection{The linear case}
We consider two-dimensional linear KP-type equations
\begin{equation}\label{2.1}
(u_t-Lu_x)_x+u_{yy}=0,\quad u(0,x,y)=\varphi(x,y)\, ,
\end{equation}
where 
\begin{equation}\label{2.2}
\widehat{Lf}(\xi)=\varepsilon|\xi|^{\alpha}\widehat{f}(\xi),\quad \xi\in \R,
\end{equation}
where $\varepsilon=1$ (KP-II type equations) or $\varepsilon=-1$ (KP-I type
equations). We denote by $G$ the fundamental solution
$$
G(t,x,y)={\mathcal F}^{-1}_{(\xi,\eta)\rightarrow (x,y)}
\big[
e^{it(\varepsilon\xi|\xi|^{\alpha}-\eta^2/\xi)}
\big]
$$
A priori, we only have that $G(t,\cdot,\cdot)\in {\mathcal S}'(\R^2)$.
Actually, for $t\neq 0$, $G(t,\cdot,\cdot)$ has a very particular form which
is the main result of this section.
\begin{theo}\label{thm2.1}
Suppose that $\alpha>1/2$ in (\ref{2.2}).
Then for $t\neq 0$,
$$
G(t,\cdot,\cdot)\in  C(\R^2)\cap L^{\infty}(\R^2).
$$
Moreover, for $t\neq 0$, there exists 
$$
A(t,\cdot,\cdot)\in C(\R^2)\cap L^{\infty}(\R^2)\cap C^{1}_{x}(\R^2)
$$
($C^{1}_{x}(\R^2)$ denotes the set of $x$ differentiable continuous
function on $\R^2$) such that
$$
G(t,x,y)=\frac{\partial A}{\partial x}(t,x,y).
$$
In addition for $t\neq 0$, $y\in\R$, $\varphi\in L^1(\R^2)$,
$$
\lim_{|x|\rightarrow \infty}(A\star\varphi)(t,x,y)=0.
$$
As a consequence, the solution of (\ref{2.1}) with data $\varphi\in L^1(\R^2)$
is given by 
$$
u(t,\cdot,\cdot)\equiv S(t)\varphi=G\star \varphi
$$
and
$$
u(t,\cdot,\cdot)=\frac{\partial}{\partial x}\big( A\star \varphi\big).
$$
One has therefore
$$
\int_{-\infty}^{\infty}u(t,x,y)\, dx =0,\quad \forall y\in\R,\,\,\, \forall t\neq
0\, ,
$$
in the sense of generalized Riemann integrals.
\end{theo}
\begin{rem}
It is worth noticing that the result of Theorem~\ref{thm2.1} is related to the
infinite speed of propagation of the KP free evolutions.
Let us also notice that the assumption $\alpha>1/2$ can be relaxed, if we
assume that a sufficient number of derivatives of $\varphi$ belong to $L^1$.
Such an assumption is however not natural in the context of the KP equations.
\end{rem}
\begin{rem}
In the case of the classical KP-II equation ($\alpha=2$, $\varepsilon=+1$),
Theorem~\ref{thm2.1} follows from an observation of Redekopp
\cite{R}. Namely, one has
$$
G(t,x,y)=-\frac{1}{3t}{\rm Ai}(\zeta)\, {\rm Ai}'(\zeta),
$$
where ${\rm Ai}$ is the Airy function and
$$
\zeta=c_1\frac{x}{t^{1/3}}+c_2\frac{y^2}{t^{4/3}}
$$
for some real constants $c_1> 0$ and $c_2>0$.
Thus 
$G(t,x,y)=\frac{\partial}{\partial x}A(t,x,y)$
with
$$
A(t,x,y)=-\frac{1}{6c_1 t^{2/3}}\, {\rm Ai}^{2}\Big( 
c_1\frac{x}{t^{1/3}}+c_2\frac{y^2}{t^{4/3}}\Big)
$$
and
$$
u=\frac{\partial A}{\partial x}\star \varphi=\frac{\partial}{\partial x}\big(
A\star \varphi\big),
$$
which proves the claim for the KP-II equation
(the fact that $\lim_{|x|\rightarrow \infty}A(t,x,y)=0$ results from 
a well known decay property of the  Airy function).
A similar explicit computation does not seem to be valid for the classical
KP-I equation or for KP type equations with general symbols.
\end{rem}
\begin{proof}[Proof of Theorem~\ref{thm2.1}.]
We will consider only the case $\varepsilon=1$ in (\ref{2.1}). The analysis in the case $\varepsilon=-1$
is analogous. It is plainly sufficient to consider only the case $t>0$. We have
\begin{equation}\label{pak}
G(t,x,y)=
(2\pi)^{-2}\int_{\R^2}
e^{i(x\xi+y\eta)+it(\xi|\xi|^{\alpha}-\eta^2/\xi)}d\xi d\eta\, ,
\end{equation}
where the last integral has the usual interpretation.
We first check that $G(t,x,y)$ is a continuous function of $x$
and $y$. By the change of variables 
$$
\eta'=\frac{t^{1/2}}{|\xi|^{1/2}}\eta,
$$
we obtain
\begin{multline*}
G(t,x,y) =  \frac{c}{t^{1/2}}
\int_{\R_{\xi}}
|\xi|^{1/2}
\Big(
\int_{\R_{\eta}}
e^{i(y/t^{1/2})|\xi|^{1/2}\eta-i{\rm sgn}(\xi)\eta^2}d\eta
\Big)
e^{ix\xi +it\xi|\xi|^{\alpha}} d\xi
\\
 = 
\frac{c}{t^{1/2}}
\int_{\R}
e^{-i({\rm sgn}(\xi))\frac{\pi}{4}}\,  
|\xi|^{1/2}e^{iy^2\xi/4t}\,
e^{ix\xi +it\xi|\xi|^{\alpha}} d\xi\, 
\\ 
= 
\frac{c}{t^{1/2}}
\int_{\R}
e^{-i({\rm sgn}(\xi))\frac{\pi}{4}}\,  
\,  |\xi|^{1/2}
e^{i\xi(x+ y^2/4t)}\, e^{it\xi |\xi|^{\alpha}} d\xi
\\
=
\frac{c}{t^{\frac{1}{2}+\frac{3}{2(\alpha+1)}}}
\int_{\R}
e^{-i({\rm sgn}(\xi))\frac{\pi}{4}}\,  
|\xi|^{1/2}
\exp
\Big(
i\xi
\big(
\frac{x}{t^{\frac{1}{\alpha+1}}}+\frac{y^2}{4t^{\frac{\alpha+2}{\alpha+1}}}
\big)
\Big)
e^{i\xi |\xi|^{\alpha}} d\xi\,.
\end{multline*}
Let us define
$$
H(\lambda)=c
\int_{\R}
e^{-i({\rm sgn}(\xi))\frac{\pi}{4}}\,  |\xi|^{1/2}
e^{i\lambda \xi}\, e^{i\xi |\xi|^{\alpha}} d\xi\, .
$$
Then $H$ is continuous in $\lambda$. We will only consider the worse case
$\lambda\leq 0$.
The phase $\varphi(\xi)=i(\lambda\xi+\xi |\xi|^{\alpha})$ has then $2$
critical points $\pm\xi_{\alpha}$ where
$\xi_{\alpha}=\Big(\frac{\mu}{\alpha+1}\Big)^{1/\alpha}$,
$\mu=-\lambda$. We write, for $\varepsilon>0$ small enough
$$
H(\lambda)=
\int_{-\infty}^{-\xi_{\alpha}-\varepsilon}\cdots
+\int_{-\xi_{\alpha}-\varepsilon}^{\xi_{\alpha}+\varepsilon}\cdots
+
\int_{\xi_{\alpha}+\varepsilon}^{\infty}\cdots:=I_1(\lambda)+I_2(\lambda)+I_3(\lambda)\, .
$$
Clearly $I_2(\lambda)$ is a  continuous function of $\lambda$. We consider only $I_3(\lambda)$,
\begin{multline*}
I_3(\lambda)=
c\int_{\xi_{\alpha}+\varepsilon}^{\infty}\frac{\xi^{1/2}}{\varphi'(\xi)}\,
\frac{d}{d\xi}\Big[e^{\varphi(\xi)}\Big]\,d\xi
\\
=
c\Big[
\frac{\xi^{1/2}e^{\varphi(\xi)}}
{\lambda+\xi^{\alpha}(\alpha+1)}
\Big]^{\infty}_{\xi_{\alpha}+\varepsilon}
+c\int_{\xi_{\alpha}+\varepsilon}^{\infty}
\Big[
\frac{1}{2(\lambda+\xi^{\alpha}(\alpha+1))\xi^{1/2}}
-
\frac{\alpha(\alpha+1)\xi^{\alpha-1/2}}
{(\lambda+(\alpha+1)\xi^{\alpha})^{2}}
\Big]e^{\varphi(\xi)}d\xi
\end{multline*}
which for $\alpha>1/2$ defines a continuous function of $\lambda$.
Hence the integral (\ref{pak}) is a continuous function of $(x,y)$ which
coincides with the inverse Fourier transform (in ${\mathcal S}'(\R^2)$) of $\exp(it(\xi|\xi|^{\alpha}-\eta^2/\xi))$.
\\

We next set for $t>0$,
$$
A(t,x,y)\equiv(2\pi)^{-2}
\int_{\R^2}
\frac{1}{i\xi}e^{i(x\xi+y\eta)+it(\xi|\xi|^{\alpha}-\eta^2/\xi)}d\xi d\eta\, .
$$
The last integral is clearly not absolutely convergent not only at infinity
but also for $\xi$ near zero. Nevertheless, the oscillations involved in its
definition will allow us to show that $A(t,x,y)$ is in fact a
continuous function. By the change of variables 
$$
\eta'=\frac{t^{1/2}}{|\xi|^{1/2}}\eta,
$$
we obtain
\begin{multline*}
A(t,x,y) =  \frac{c}{t^{1/2}}
\int_{\R_{\xi}}
\frac{{\rm sgn}(\xi)}{|\xi|^{1/2}}
\Big(
\int_{\R_{\eta}}
e^{i(y/t^{1/2})|\xi|^{1/2}\eta-i{\rm sgn}(\xi)\eta^2}d\eta
\Big)
e^{ix\xi +it\xi|\xi|^{\alpha}} d\xi
\\
 = 
\frac{c}{t^{1/2}}
\int_{\R}
\frac{
({\rm sgn}(\xi))
e^{-i({\rm sgn}(\xi))\frac{\pi}{4}}
}{|\xi|^{1/2}}
e^{iy^2\xi/4t}\,
e^{ix\xi +it\xi|\xi|^{\alpha}} d\xi\, 
\\ 
= 
\frac{c}{t^{\frac{\alpha+2}{2(\alpha+1)}}}
\int_{\R}
\frac{
({\rm sgn}(\xi))
e^{-i({\rm sgn}(\xi))\frac{\pi}{4}}
}{|\xi|^{1/2}}
\exp
\Big(
i\xi
\big(
\frac{x}{t^{\frac{1}{\alpha+1}}}+\frac{y^2}{4t^{\frac{\alpha+2}{\alpha+1}}}
\big)
\Big)
e^{i\xi|\xi|^{\alpha}}d\xi\,.
\end{multline*}
We now need the following lemma.
\begin{lem}\label{vaj}
Let for $\alpha>0$
$$
F(\lambda)=\int_{\R}
\frac{({\rm sgn}(\xi))e^{-i({\rm sgn}(\xi))\frac{\pi}{4}}
}{|\xi|^{1/2}}
e^{i\lambda\xi+i\xi |\xi|^{\alpha}}d\xi\, .
$$
Then $F$ is a continuous function which tends to zero as $|\lambda|\rightarrow
+\infty$. 
\end{lem}
\begin{proof}
Write $F$ as
$$
F(\lambda)=\int_{|\xi|\leq 1}\dots +
\int_{|\xi|\geq 1}\dots := F_{1}(\lambda)+F_{2}(\lambda)\, .
$$
Since $|\xi|^{-1/2}$ in integrable near the origin, by the Riemann-Lebesgue lemma
$F_{1}(\lambda)$ is continuous and 
$$
\lim_{|\lambda|\rightarrow \infty} F_{1}(\lambda)=0\, .
$$
We consider two cases in the analysis of $F_{2}(\lambda)$
\\

{\bf Case 1.} $\lambda\geq -1$.
\\
After an integration by parts, we obtain that
\begin{multline}\label{inte}
F_{2}(\lambda)=
\frac{c\cos(\lambda+1-\frac{\pi}{4})}{\lambda+\alpha+1}+
\\
+c\int_{1}^{\infty}
\cos\Big(\lambda\xi+\xi^{\alpha+1}-\frac{\pi}{4}\Big)\,\,
\frac{\lambda+(\alpha+1)(2\alpha+1)\xi^{\alpha}}
{\xi^{3/2}(\lambda+(\alpha+1)\xi^{\alpha})^{2}}\,d\xi\,.
\end{multline}
The first term is clearly a continuous function of $\lambda$ which tends to
zero as $\lambda\rightarrow\infty$. 
Observing that
$$
0\leq 
\frac{\lambda+(\alpha+1)(2\alpha+1)\xi^{\alpha}}
{\xi^{3/2}(\lambda+(\alpha+1)\xi^{\alpha})^{2}}
\leq C_{\alpha}\xi^{-3/2}
$$
uniformly with respect to $\xi\geq 1$ and $\lambda\geq -1$, we deduce from 
the dominated convergence theorem that the right hand-side of (\ref{inte}) is a
continuous function of $\lambda$ for $\lambda\geq -1$.
On the other hand for $\lambda\geq 1$,
$$
\frac{\lambda+(\alpha+1)(2\alpha+1)\xi^{\alpha}}
{\xi^{3/2}(\lambda+(\alpha+1)\xi^{\alpha})^{2}}
\leq
\frac{2\alpha+1}{\lambda\xi^{3/2}}
$$
and thus the right hand-side of (\ref{inte}) tends to zero as $\lambda\rightarrow +\infty$.
\\

{\bf Case 2.} $\lambda\leq -1$.
\\ 
Set $\lambda=-\mu$ with $\mu\geq 1$. In the integral over $|\xi|\geq 1$ defining
$F_{2}(\lambda)$,
we consider only the integration over $[1,+\infty[$. The integration over $]-\infty,-1]$ 
can be treated in a completely analogous way. We perform the changes of
variables
$$
\xi\longrightarrow \xi^2
$$
and
$$
\xi\longrightarrow \mu^{\frac{1}{2\alpha}}\xi
$$
to conclude that
$$
\widetilde{F_2}(\lambda):=
c\int_{1}^{\infty}\frac{1}{\xi^{1/2}}e^{i\lambda\xi+i\xi |\xi|^{\alpha}}d\xi
=
c\mu^{\frac{1}{2\alpha}}
\int_{\mu^{-\frac{1}{2\alpha}}}^{\infty}
e^{i\mu^{1+\frac{1}{\alpha}}\big[\xi^{2(\alpha+1)}-\xi^{2}\big]}
d\xi
$$
Let us set 
$$
\varphi(\xi)=\xi^{2(\alpha+1)}-\xi^{2}\, .
$$
Then
$$
\varphi'(\xi)=2\xi[(\alpha+1)\xi^{2\alpha}-1]\,.
$$
Let us split
$$
\widetilde{F_2}(\lambda)=
c\mu^{\frac{1}{2\alpha}}
\int_{\mu^{-\frac{1}{2\alpha}}}^{1}\dots+
c\mu^{\frac{1}{2\alpha}}\int_{1}^{\infty}\dots
:=
I_1(\mu)+I_2(\mu)\, .
$$
Since $\varphi'(\xi)$ does not vanish for $\xi\geq 1$, we can integrate by parts which gives
$$
I_2(\mu)=
\frac{1}{2i\mu^{1+\frac{1}{2\alpha}}}
\Big(
\frac{c}{\alpha}+
c\int_{1}^{\infty}
e^{i\mu^{1+\frac{1}{\alpha}}\big[\xi^{2(\alpha+1)}-\xi^{2}\big]}\,\,
\frac{(\alpha+1)(2\alpha+1)\xi^{2\alpha}-1}
{\xi^{2}((\alpha+1)\xi^{2\alpha}-1)^{2}}
d\xi
\Big)
$$
which is a continuous function of $\mu\geq 1$ thanks to the 
dominated convergence theorem. Moreover it clearly tends to zero as
$\mu\rightarrow+\infty$.
\\

Let us next analyse $I_{1}(\mu)$. We first observe
that thanks to the dominated convergence theorem, $I_{1}(\mu)$ is a continuous
function of $\mu$. It remains to prove that 
$I_{1}(\mu)\rightarrow 0$ as $\mu\rightarrow \infty$.
For $\xi\in [\mu^{-\frac{1}{2\alpha}},1]$, the phase $\varphi$ has a critical
point and a slightly more delicate argument is needed. Compute
$$
\varphi''(\xi)=2[(\alpha+1)(2\alpha+1)\xi^{2\alpha}-1]\, .
$$
Observe that $\varphi'(\xi)$ is vanishing only at zero and
$$
\xi_1(\alpha):=\left(\frac{1}{\alpha+1}\right)^{\frac{1}{2\alpha}}\,\, .
$$
Next, we notice that $\varphi''(\xi)$ is vanishing at
$$
\xi_2(\alpha):=\left(\frac{1}{(\alpha+1)(2\alpha+1)}\right)^{\frac{1}{2\alpha}}\,\, .
$$
Clearly $\xi_2(\alpha)<\xi_1(\alpha)<1$ and we choose a real number $\delta$
such that
$$
\xi_2(\alpha)<\delta<\xi_1(\alpha)<1\, .
$$
For $\mu\gg 1$, we can split
$$
I_{1}(\mu)=c\mu^{\frac{1}{2\alpha}}\int_{\mu^{-\frac{1}{2\alpha}}}^{\delta}\dots
+
c\mu^{\frac{1}{2\alpha}}\int_{\delta}^{1}
:= J_{1}(\mu)+J_{2}(\mu)\,\, .
$$
For $\xi\in [\mu^{-\frac{1}{2\alpha}},\delta]$, we have the minoration
$$
|\varphi'(\xi)|\geq c\mu^{-\frac{1}{2\alpha}}>0
$$
and an integration by parts shows that
$$
J_{1}(\mu)=\mu^{\frac{1}{\alpha}}{\mathcal O}(\mu^{-1-\frac{1}{\alpha}})
\leq C\mu^{-1}
$$
which clearly tends to zero as $\mu\rightarrow \infty$.
For $\xi\in [\delta,1]$, we have the 
minoration
$$
|\varphi''(\xi)|\geq c>0
$$
and therefore we can apply 
the Van der Corput lemma (see \cite[ Proposition~2]{stein}) 
to conclude that
$$
J_{2}(\mu)=\mu^{\frac{1}{2\alpha}}{\mathcal O}(\mu^{-\frac{1}{2}-\frac{1}{2\alpha}})
\leq C\mu^{-\frac{1}{2}}
$$
which tends to zero as $\mu\rightarrow \infty$.
This completes the proof of Lemma~\ref{vaj}.
\end{proof}
It is now easy to check that $\partial_{x}A=G$ in the sense of distributions. Since both $A$ and $G$ are
continuous, we deduce that $A$ has a classical derivative with
respect to $x$ which is equal to $G$. Finally,
since $\varphi\in L^1(\R^2)$, applying Lemma ~\ref{vaj} and the Lebesgue Theorem completes the proof of Theorem~\ref{thm2.1}.
\end{proof}
%%%%%%%%%%%%%%%%%%%%%%%%%%%%%%%%%%%%%%%%%%%%%%%%%%%%%%%%%%%%%%%%%%%%%%%%%%%%%%%%%%%%%%%%%%%%%%%%%%%%%%%%%%%%%
%%%%%%%%%%%%%%%%%%%%%%%%%%%%%%%%%%%%%%%%%%%%%%%%%%%%%%%%%%%%%%%%%%%%%%%%%%%%%%%%%%%%%%%%%%%%%%%%%%%%%%%%%%%%
\subsection{The nonlinear case}
After a change of frame we can eliminate the $u_x$ term and reduce the Cauchy problem for (\ref{1.9}) to
\begin{equation}\label{2.3}
(u_t+uu_x-Lu_x)_x+u_{yy}=0,\quad u(0,x,y)=\varphi(x,y)\, ,
\end{equation}
In order to state our result concerning (\ref{2.3}), for $k\in \N$, we denote by $H^{k,0}(\R^2)$ the Sobolev space
of $L^2(\R^2)$ functions $u(x,y)$ such that $\partial_{x}^{k}u\in L^2(\R^2)$.
\begin{theo}\label{thm2.2}
Assume that $\alpha>1/2$. 
Let $\varphi\in L^{1}(\R^2)\cap H^{2,0}(\R^2)$ and 
\begin{equation}\label{new}
u\in C([0,T]\, ;\, H^{2,0}(\R^2))
\end{equation}
be a distributional solution of (\ref{2.3}).
Then, for every $t\in (0,T]$, $u(t,\cdot,\cdot)$ is a continuous function of
$x$ and $y$ which satisfies
\begin{equation*}
%\label{2.4}
\int_{-\infty}^{\infty}u(t,x,y)dx =0,\quad \forall y\in \R,\,\,\, \forall t\in (0,T]
\end{equation*}
in the sense of generalized Riemann integrals.
Moreover, $u(t,x,y)$ is the
derivative with respect to $x$ of a $C^1_x$ continuous function
which vanishes as
$x\rightarrow \pm \infty$ for every fixed $y\in\R$ and $t\in [0,T]$.
\end{theo}
\begin{rem}
The case $\alpha=2$ corresponds to the classical KP-I, KP-II equations.
In the case of the KP-II, we have global solutions for data in $L^{1}(\R^2)\cap H^{2,0}(\R^2)$
(see \cite{B}). Thus Theorem~\ref{thm2.2} displays a striking smoothing
effect of the KP-II equation : $u(t,\cdot,\cdot)$ becomes a continuous
function of $x$ and $y$ 
(with zero mean in $x$) 
for $t\neq 0$
(note that $L^{1}(\R^2)\cap H^{2,0}(\R^2)$ in not included in $C^{0}(\R^2)$).
A similar comment is valid for the local solutions of the KP-I equation
in \cite{MST1} and more especially in \cite{MST3}.
\end{rem}
\begin{rem}
The numerical simulations in \cite{KSM} display clearly the phenomena described in Theorem
\ref{thm2.2} in the case of the KP-I equation.
\end{rem}
\begin{proof}[Proof of Theorem~\ref{thm2.2}]
Under our assumption on $u$, one has the Duhamel representation
\begin{equation}\label{2.5}
u(t)=S(t)\varphi-\int_{0}^{t}S(t-s)\big(u(s)u_{x}(s)\big)ds
\end{equation}
where
\begin{multline*}
%\label{2.5pak}
\int_{0}^{t}S(t-s)\big(u(s)u_{x}(s)\big)ds
=
\\
=
\int_{0}^{t}\partial_{x}
\Big(
\int_{\R^2}A(x-x',y-y',t-s)(uu_x)(x',y',s)dx'dy'
\Big)ds\,.
\end{multline*}
>From Theorem~\ref{thm2.1}, it suffices to consider only the integral
term in the right hand-side of (\ref{2.5}). Using the notations of
Lemma~\ref{vaj},
$$
A(x-x',y-y',t-s)=\frac{c}{(t-s)^{\frac{\alpha+2}{2(\alpha+1)}}}
F\Big(
\frac{x-x'}{(t-s)^{\frac{1}{\alpha+1}}}+\frac{(y-y')^2}{4(t-s)^{\frac{\alpha+2}{\alpha+1}}}
\Big).
$$
Recall that $F$ is a continuous and bounded function on $\R$. 
Next we set
$$
I(x,y,t-s,s)\equiv\partial_{x}
\Big(
\int_{\R^2}A(x-x',y-y',t-s)(uu_x)(x',y',s)dx'dy'
\Big).
$$
Using the Lebesgue differentiation theorem and the assumption (\ref{new}), we can write
$$
I(x,y,t-s,s)=
\int_{\R^2}A(x-x',y-y',t-s)\partial_{x}(uu_x)(x',y',s)dx'dy' \,.
$$
Moreover for $\alpha>0$,
$$
\frac{\alpha+2}{2(\alpha+1)}<1
$$
and therefore $I$ is integrable in $s$ on $[0,t]$.
Therefore, by the Lebesgue differentiation theorem,
\begin{equation}\label{function}
\int_{0}^{t}
\int_{\R^2}A(x-x',y-y',t-s)(uu_x)(x',y',s)dx'dy'
ds
\end{equation}
is a $C^{1}_{x}$ function and
\begin{multline*}
\int_{0}^{t}S(t-s)\big(u(s)u_{x}(s)\big)ds
=
\\
=
\partial_{x}\Big(
\int_{0}^{t}
\int_{\R^2}A(x-x',y-y',t-s)(uu_x)(x',y',s)dx'dy'
ds\Big)
\end{multline*}
Let us finally show that for fixed $y$ and $t$ the function (\ref{function})
tends to zero as $x$ tends to $\pm\infty$. 
For that purpose it suffices to apply the Lebesgue dominated convergence
theorem to the integral in $s$, $x'$, $y'$. Indeed, for fixed $s$, $x'$, $y'$,
the function under the integral tends to  zero as $x$ tends to $\pm\infty$
thanks to the linear analysis. 
On the other hand, using Lemma~\ref{vaj}, we can write
$$
|A(x-x',y-y',t-s)(uu_x)(x',y',s)|\leq 
\frac{c}{(t-s)^{\frac{\alpha+2}{2(\alpha+1)}}}
|(uu_x)(x',y',s)|\, .
$$
Thanks to the assumptions on $u$ the right hand-side of the above inequality
is integrable in $s$, $x'$, $y'$ and independent of $x$. 
Thus we can apply the Lebesgue dominated convergence
theorem to conclude that the function (\ref{function})
tends to zero as $x$ tends to $\pm\infty$. 
This completes the proof of Theorem~\ref{thm2.2}.
\end{proof}
\begin{remarque}
If $\alpha>2$, the assumptions can be weakened to $\varphi\in L^{1}(\R^2)\cap
H^{1,0}(\R^2)$ and $u\in C([0,T];H^{1,0}(\R^2))$.
This results follows from the fact that the fundamental solution $G$
write
$$
G(t,x,y)=
\frac{c}{t^{1/2+3/2(\alpha+1)}}
B(t,x,y),
$$
where $B\in L^{\infty}$ and $1/2+3/(2(\alpha+1)<1$ for $\alpha>2$.
\end{remarque}
%%%%%%%%%%%%%%%%%%%%%%%%%%%%%%%%%%%%%%%%%%%%%%%%%%%%%%%%%%%%%%%%%%%%%%%%%%%%%%%%%%%%%%%%%%%%%%%%%%%%%%%%%%%%%%%%%%%%%%%%
%%%%%%%%%%%%%%%%%%%%%%%%%%%%%%%%%%%%%%%%%%%%%%%%%%%%%%%%%%%%%%%%%%%%%%%%%%%%%%%%%%%%%%%%%%%%%%%%%%%%%%%%%%%%%%%%%%%%%%%%
\section{KP-BBM type equations}
\subsection{The linear case}
We consider the Cauchy problem
\begin{equation}\label{3.1}
(u_t+u_x+Lu_t)_x+u_{yy}=0,\quad u(0,x,y)=\varphi(x,y)\, ,
\end{equation}
where $L$ is given by (\ref{2.2}) with $\varepsilon=1$. 
A simple computation shows that the fundamental solution is given by
$$
G(t,x,y)={\mathcal
  F}^{-1}\big[e^{-i\frac{t}{1+|\xi|^{\alpha}}(\xi+\eta^2/\xi)}\big]\, .
$$
Due to the bad oscillatory properties of the phase we have to modify a little
bit the statement of Theorem~\ref{thm2.1}.
\begin{theo}\label{thm3.1}
Assume that $\alpha>0$. Let $\varphi$ be such that
$(I-\partial_{x}^{2})^{\beta/2}\varphi\in L^1(\R^2)$
with $\beta>(\alpha+3)/2$. Then 
$$
u(t,\cdot,\cdot)\equiv S(t)\varphi=G\star \varphi
$$
can be written as
$$
\frac{\partial}{\partial x}\big( A\star (I-\partial_{x}^{2})^{\beta/2}\varphi\big),
$$
where
$$
A=\partial_{x}^{-1}(I-\partial_{x}^{2})^{-\beta/2}G
$$
and
$$
A(t,\cdot,\cdot)\in C(\R^2)\cap L^{\infty}(\R^2)\cap C^1_{x}(\R^2)  \, .
$$
Moreover, $(A\star  (I-\partial_{x}^{2})^{\beta/2}\varphi)(t,x,y)$ is, for fixed $t\neq 0$ continuous in $x$
and $y$ and satisfies
$$
\lim_{|x|\rightarrow \infty}(A\star (I-\partial_{x}^{2})^{\beta/2}\varphi)(t,x,y)=0
$$
for any $y\in\R$ and $t\neq 0$.
Thus,
$$
\int_{-\infty}^{\infty}u(t,x,y)\, dx =0,\quad \forall y\in\R,\,\,\, \forall t\neq
0\, ,
$$
in the sense of generalized Riemann integrals.
\end{theo}
\begin{rem}
The assumption $(I-\partial_{x}^{2})^{\beta/2}\varphi\in L^1(\R^2)$, is
natural in the context of KP-BBM problems, in view of the weak dispersive
properties of the free BBM evolution.
\end{rem}
\begin{proof}
We set
$$
\tilde{G}(t,x,y)=\int_{\R^2}
\frac{1}{(1+|\xi|^2)^{\beta/2}}\, \,
e^{-\frac{it}{1+|\xi|^{\alpha}}(\xi+\eta^2/\xi)}
e^{ix\xi+iy\eta}d\xi d\eta\, .
$$
Setting
$$
\eta'=\frac{t^{1/2}\eta}{|\xi|^{1/2}(1+|\xi|^{\alpha})^{1/2}}\, ,
$$
we obtain
\begin{multline*}
\tilde{G}(t,x,y)=
\\
\frac{c}{t^{1/2}}\int_{\R_{\xi}}
|\xi|^{1/2}
\frac{(1+|\xi|^{\alpha})^{1/2}}{(1+|\xi|^{2})^{\beta/2}}\Big(\int_{\R_{\eta}}
e^{-i{\rm sgn}(\xi)\eta^2}
e^{i(y\eta/t^{1/2})|\xi|^{1/2}(1+|\xi|^{\alpha})^{1/2}}
d\eta
\Big)
e^{ix\xi-i\frac{t\xi}{1+|\xi|^{\alpha}}}
d\xi
\\
=\frac{c}{t^{1/2}}
\int_{\R}
|\xi|^{1/2}
e^{-i({\rm sgn}(\xi))\frac{\pi}{4}}
\frac{(1+|\xi|^{\alpha})^{1/2}}{(1+|\xi|^{2})^{\beta/2}}
e^{-i\frac{t\xi}{1+|\xi|^{\alpha}}}
e^{i\lambda\xi}\, e^{i(y^{2}/4t) \xi|\xi|^{\alpha}}d\xi\, ,
\end{multline*}
where $\lambda=x+y^2/4t$.
Setting
$$
a_{t}(\xi)=
|\xi|^{1/2}
e^{-i({\rm sgn}(\xi))\frac{\pi}{4}}
\frac{(1+|\xi|^{\alpha})^{1/2}}{(1+|\xi|^{2})^{\beta/2}}
e^{-it\frac{\xi}{1+|\xi|^{\alpha}}}
$$
we clearly have
$$
\tilde{G}(t,x,y)
=
\frac{c}{t^{1/2}}
\int a_{t}(\xi) e^{i\lambda\xi}\, e^{i(y^{2}/4t)\xi|\xi|^{\alpha}}d\xi\, .
$$
We have the following lemma.
\begin{lem}
Let us fix $y\in\R$ and $t>0$. Set
$$
F(\lambda):=\int a_{t}(\xi) e^{i\lambda\xi}\, e^{i(y^2/4t)\xi|\xi|^{\alpha}}d\xi\,.
$$
Then $F$ is a continuous function such that
$$
\lim_{|\lambda|\rightarrow \infty}F(\lambda)=0\, .
$$
\end{lem}
\begin{proof}
It suffices to apply the Riemann-Lebesgue lemma since $a\in L^1(\R)$ when
$\beta>(\alpha+3)/2$.
\end{proof}
Next, we set
$$
\tilde{A}(t,x,y)=-i\int_{\R^2}
\frac{1}{\xi(1+|\xi|^2)^{\beta/2}}\, \,
e^{-\frac{it}{1+|\xi|^{\alpha}}(\xi+\eta^2/\xi)}
e^{ix\xi+iy\eta}d\xi d\eta\, .
$$
Similarly to above, we set
$$
\eta'=\frac{\eta t^{1/2}}{|\xi|^{1/2}(1+|\xi|^{\alpha})^{1/2}}\, ,
$$
and therefore
\begin{multline*}
\tilde{A}(t,x,y)
=\frac{c}{t^{1/2}}\int_{\R_{\xi}}
\frac{{\rm sgn}(\xi)}{|\xi|^{1/2}}
\frac{(1+|\xi|^{\alpha})^{1/2}}{(1+|\xi|^{2})^{\beta/2}}
\Big(
\int_{\R_{\eta}}
e^{-i{\rm sgn}(\xi)\eta^2}
e^{iy\eta/t^{1/2}|\xi|^{1/2}(1+|\xi|^{\alpha})^{1/2}}
d\eta
\Big)
e^{ix\xi-i\frac{t\xi}{1+|\xi|^{\alpha}}}
d\xi
\\
=\frac{c}{t^{1/2}}
\int_{\R}
\frac{{\rm sgn}(\xi)}{|\xi|^{1/2}}
e^{-i({\rm sgn}(\xi))\frac{\pi}{4}}
\frac{(1+|\xi|^{\alpha})^{1/2}}{(1+|\xi|^{2})^{\beta/2}}
e^{-it\frac{\xi}{1+|\xi|^{\alpha}}}
e^{i\lambda\xi}\, e^{i(y^{2}/4t) \xi|\xi|^{\alpha}}d\xi\, ,
\end{multline*}
where $\lambda=x+y^2/4t$. Setting
$$
\tilde{a}_t(\xi)=\frac{{\rm sgn}(\xi)}{|\xi|^{1/2}}
e^{-i({\rm sgn}(\xi))\frac{\pi}{4}}
\frac{(1+|\xi|^{\alpha})^{1/2}}{(1+|\xi|^{2})^{\beta/2}}
e^{-i\frac{t\xi}{1+|\xi|^{\alpha}}}
$$
we clearly have
$$
\tilde{A}(t,x,y)
=
\frac{c}{t^{1/2}}
\int \tilde{a}_t(\xi) e^{i\lambda\xi}\, e^{i(y^{2}/4t)\xi|\xi|^{\alpha}}d\xi\, .
$$
\begin{lem}\label{lemme3.1}
Let us fix $y\in\R$ and $t>0$. Set
$$
F_1(\lambda):=\int \tilde{a}_t(\xi) e^{i\lambda\xi}\, e^{i(y^2/4t)\xi|\xi|^{\alpha}}d\xi\, .
$$
Then $F_1$ is a continuous function such that
$$
\lim_{|\lambda|\rightarrow \infty}F_1(\lambda)=0\, .
$$
\end{lem}
\begin{proof}
It suffices to apply the Riemann-Lebesgue lemma since $\tilde{a}_t\in L^1(\R)$.
\end{proof}
The proof of Theorem~\ref{thm3.1} is now straightforward.
\end{proof}
%%%%%%%%%%%%%%%%%%%%%%%%%%%%%%%%%%%%%%%%%%%%%%%%%%%%%%%%%%%%%%%%%%%%%%%%%%%%%%%%%%%%%%%%%%%%%%%%%%%%%%%%%%
%%%%%%%%%%%%%%%%%%%%%%%%%%%%%%%%%%%%%%%%%%%%%%%%%%%%%%%%%%%%%%%%%%%%%%%%%%%%%%%%%%%%%%%%%%%%%%%%%%%%%%%%%%
%%%%%%%%%%%%%%%%%%%%%%%%%%%%%%%%%%%%%%%%%%%%%%%%%%%%%%%%%%%%%%%%%%%%%%%%%%%%%%%%%%%%%%%%%%%%%%%%%%%%%%%%%%
\subsection{The nonlinear case}
We investigate the Cauchy problem
\begin{equation}\label{3.2}
(u_t+u_x+uu_x+Lu_t)_x+u_{yy}=0,\quad u(0,x,y)=\varphi(x,y)\, .
\end{equation}
\begin{theo}\label{thm3.2}
Let $\alpha>0$, $k>\frac{\alpha+3}{4}$.
Assume that $(I-\partial_{x}^{2})^{k}\varphi\in L^1(\R^2)\cap L^2(\R^2)$.
Let $u$ be a solution of (\ref{3.2}) such that
$$
u\in C([0,T]\, ;\, H^{2k+1,0}(\R^2))\, .
$$
Then, for any $t\in (0,T]$, $u(t,x,y)$ is a continuous function in $x$ and
$y$ and satisfies
\begin{equation*}
%\label{3.3}
\int_{-\infty}^{\infty}u(t,x,y)dx =0,\quad \forall y\in \R,\,\,\, \forall t\in (0,T]
\end{equation*}
in  the sense of generalized Riemann integrals. In fact $u(t,x,y)$ is the
derivative with respect to $x$ of a $C^1_{x}$ continuous function which vanishes as
$x\rightarrow \pm \infty$, for any fixed $y$ and $t\in(0,T]$.
\end{theo}
\begin{proof}
Again we use the Duhamel formula
$$
u(t)=S(t)\varphi-\int_{0}^{t}S(t-s)u(s)u_{x}(s)ds,
$$
where
$$
S(t)=e^{it(I+L)^{-1}(\partial_x+\partial_{x}^{-1}\partial_{y}^{2})}\, .
$$
By Theorem~\ref{thm3.1}, it suffices to consider the integral term in the
Duhamel formula. 
We have
$$
\int_{0}^{t}S(t-s)u(s)u_{x}(s)ds=
\frac{1}{2}
\int_{0}^{t}\frac{\partial}{\partial x}\Big(S(t-s)u^{2}(s)\Big)ds
=
\frac{1}{2}
\frac{\partial}{\partial x}
\int_{0}^{t}S(t-s)u^{2}(s)ds.
$$
To justify the last equality, we have to check that $S(t-s)u(s)u_{x}(s)$ is
dominated by a $L^1(0,t)$ function uniformly in $(x,y)$. 
We write
\begin{equation}\label{3.4}
S(t-s)(u(s)u_x(s))= \tilde{A}(x,y,t-s)\star (I-\partial_x^2)^{k}u(s)u_{x}(s),
\end{equation}
where
$$
\tilde{A}(x,y,t-s)=\int_{\R^2}\frac{1}{(1+|\xi|^2)^{k}}
e^{-\frac{i(t-s)}{1+|\xi|^{\alpha}}(\xi+\eta^2/\xi)}
e^{ix\xi+iy\eta}d\xi d\eta\, .
$$
Proceeding as in the
beginning of the proof of Theorem~\ref{thm3.1}, it follows that
\begin{multline*}
\tilde{A}(x,y,t-s)=\\=
\frac{c}{(t-s)^{1/2}}
\int_{\R}\frac{
e^{-i({\rm sgn}(\xi))\frac{\pi}{4}}
|\xi|^{1/2}(1+|\xi|^{\alpha})^{1/2}}{(1+\xi^2)^k}
e^{-i\frac{\xi(t-s)}{1+|\xi|^{\alpha}}}
e^{i\lambda\xi}
e^{i(y^2/4(t-s))\xi|\xi|^{\alpha}}
d\xi\, ,
\end{multline*}
where $\lambda=x+y^2/4(t-s)$.
Since $k>\frac{\alpha+3}{4}$, the integral in $\xi$ defines a continuous
bounded function in $x,y,t,s$, by the Riemann-Lebesgue theorem. It follows
that
$$
|S(t-s)u(s)u_{x}(s)|\leq
\frac{c}{(t-s)^{1/2}}\|(I-\partial_{x}^{2})^{k}(uu_x)\|_{L^1(\R^2)}
\leq \frac{c}{(t-s)^{1/2}}
$$
since $u\in C([0,T]\, ;\, H^{2k+1,0}(\R^2)).$
Since the function
$$
\frac{e^{-i({\rm sgn}(\xi))\frac{\pi}{4}}
|\xi|^{1/2}(1+|\xi|^{\alpha})^{1/2}}{(1+\xi^2)^k}
e^{-i\frac{\xi(t-s)}{1+\xi^{\alpha}}}
$$
belongs to $L^1(\R_{\xi})$, we can use the  Riemann-Lebesgue lemma
to obtain that for fixed $y$, $t$ and $s$ the function $\tilde{A}(x,y,t-s)$
tends to zero as $x$ tends to $\pm \infty$.
Moreover the absolute value of $\tilde{A}(x,y,t-s)$ is bounded by $c|t-s|^{-1/2}$.
Thus as in the proof of Theorem~\ref{thm2.2}, we can apply the Lebesgue
dominated convergence theorem to conclude that
$$
\lim_{x\rightarrow\pm\infty}\int_{0}^{t}S(t-s)u^{2}(s)ds=0
$$
for any fixed $y\in\R$ and $t\in(0,T]$. 
This achieves the proof of Theorem~\ref{thm3.2}.
\end{proof}
\begin{rem}
For large values of $\alpha$ one can relax the assumptions on $k$ in the
hypothesis for $u$ by simply using the $H^s$ unitary property of $S(t)$.
\end{rem} 
%%%%%%%%%%%%%%%%%%%%%%%%%%%%%%%%%%%%%%%%%%%%%%%%%%%%%%%%%%%%%%%%%%%%%%%%%%%%%%%%%%%%%%%%%%%%%%%%%%%%%%%
%%%%%%%%%%%%%%%%%%%%%%%%%%%%%%%%%%%%%%%%%%%%%%%%%%%%%%%%%%%%%%%%%%%%%%%%%%%%%%%%%%%%%%%%%%%%%%%%%%%%%%%
\section{Extensions}
With the price of some technicalities, one could consider symbols $c(\xi)$ in
(\ref{1.3}) which behave like $|\xi|^\alpha$ at infinity but which are not
homogeneous.  
\\

Let us finally comment briefly on the three dimensional case. For
simplicity we consider only the KP-type equations. 
\begin{equation}\label{4.1}
(u_t-Lu_x)_{x}+u_{yy}+u_{zz}=0,\quad u(0,x,y,z)=\varphi(x,y,z)\, .
\end{equation}
with  $L$ given by (\ref{2.2}). Following the lines of the proof of
Theorem~\ref{thm2.1}, we find that the fundamental solution $G$ can be
expressed as $G=\partial_x A$, where
$$
A(t,x,y,z)=
\frac{c}{t^{1+\frac{2}{\alpha+2}}}
\int_{\R}{\rm sgn}(\xi)\, 
e^{-i({\rm sgn}(\xi))\frac{\pi}{4}}
e^{i\xi(x/t^{1/(\alpha+1)}+(y^2+z^2)/4t)}\, e^{i\xi |\xi|^{\alpha}}\,d\xi\, .
$$
We first notice that, when $\alpha>1$, $G$ is a well-defined continuous
function of $(x,y,z)$. Actually the proof follows the same lines as the $2$
dimensional case. 
Let
$$
F(\lambda)=\int_{\R}{\rm sgn}(\xi) e^{-i({\rm sgn}(\xi))\frac{\pi}{4}}e^{i\lambda \xi}\, e^{i\xi |\xi|^{\alpha}}\,d\xi\, .
$$
By a result of \cite{SSS}, $F(\lambda)$ is a continuous function which tends
to zero as $|\lambda|\rightarrow \infty$, provided $\alpha>1$
(notice that this excludes the case $\alpha=1$ which would correspond to the 3D
generalizations of the Benjamin-Ono equation).
We thus obtain the exact counterpart of Theorem~\ref{thm2.1} in the 3D case
when $\alpha>1$ (this includes the 3D usual KP equations).

\end{document}